\documentclass[12pt]{article}%
\usepackage{amsfonts}
\usepackage{graphicx}
\usepackage{amsmath}
\usepackage{amssymb}%
\providecommand{\U}[1]{\protect\rule{.1in}{.1in}}
\newtheorem{theorem}{Theorem}

\newtheorem{lemma}[theorem]{Lemma}

\begin{document}

\begin{center}
\bigskip

{\Large On the invariance principle for reversible Markov chains}

\bigskip Magda Peligrad\footnote{Department of Mathematical Sciences,
University of Cincinnati, PO Box 210025, Cincinnati, Oh 45221-0025, USA.
\texttt{ peligrm@ucmail.uc.edu}} and Sergey Utev \footnote{Department of
Mathematics, University of Leicester, University Road, LE1 7RH UK,
\texttt{\ su35@leicester.ac.uk}}
\end{center}


\bigskip

\begin{center}
Abstract
\end{center}

In this paper, we investigate the functional central limit theorem for
stochastic processes associated to partial sums of additive functionals of
reversible Markov chains with general spate space, under the normalization
standard deviation of partial sums. For this case, we show that the functional
central limit theorem is equivalent to the fact that the variance of partial
sums is regularly varying with exponent $1$ and the partial sums satisfy the
CLT. It is also equivalent to the conditional CLT.

\section{Introduction and Result}

Reversible Markov chains play a very important role in applications to
infinite particle systems, random walks, processes in random media,
Metropolis-Hastings algorithms. For instance, Kipnis and Varhadhan (1986) and
Kipnis and Landim (1999) considered applications to interacting particle
systems, Tierney (1994), Zhao et al. (2010), Longla et al. (2012) discussed
the applications to Markov Chain Monte Carlo. Our paper is motivated by the
functional limit theorem in the paper by Longla et al. (2012). Our result will
bring further clarification on this subject. Without assuming aperiodicity or
irreducibility properties, we shall show that for an additive functional of a
stationary reversible Markov chain the functional CLT is equivalent to CLT
plus the fact that the variance of partial sums is regularly varying with
exponent $1$.

We assume that $(\xi_{n})_{n\in\mathbb{Z}}$ is a stationary Markov chain
defined on a probability space $(\Omega,\mathcal{F},\mathbb{P})$ with values
in a general state space $(S,\mathcal{A})$. The marginal distribution is
denoted by $\pi(A)=\mathbb{P}(\xi_{0}\in A)$. Assume that there is a regular
conditional distribution for $\xi_{1}$ given $\xi_{0}$ denoted by
$Q(x,A)=\mathbb{P}(\xi_{1}\in A|\,\xi_{0}=x)$. Let $Q$ also denote the Markov
operator {acting via $(Qf)(x)=\int_{S}f(s)Q(x,ds).$ Next, let $\mathbb{L}%
_{0}^{2}(\pi)$ be the set of measurable functions on $S$ such that $\int
f^{2}d\pi<\infty$ and $\int fd\pi=0.$ For some function } ${f}\in$%
{$\mathbb{L}_{0}^{2}(\pi)$, let }%
\begin{equation}
{X_{i}=f(\xi_{i}),\ S_{n}=\sum\limits_{i=1}^{n}X_{i},\ }\sigma_{n}%
=(\mathbb{E}S_{n}^{2})^{1/2}. \label{defX}%
\end{equation}

For any integrable random variable $X$ we denote $\mathbb{E}_{\xi_{n}%
}(X)=\mathbb{E}(X|\xi_{n}).$ The symbol {\ }$\Rightarrow$ denotes convergence
in distribution.

The Markov chain is called reversible if $Q=Q^{\ast},$ where $Q^{\ast}$ is the
adjoint operator of \ $Q$. The condition of reversibility is equivalent to
requiring that $(\xi_{0},\xi_{1})$ and $(\xi_{1},\xi_{0})$ have the same
distribution or%
\[
\int_{A}Q(\omega,B)\pi(d\omega)=\int_{B}Q(\omega,A)\pi(d\omega),
\]
for all Borel sets $A,B\in\mathcal{A}$.

Gordin and Lif\v{s}ic (1981) proved a CLT for functionals of a normal Markov
chain. In particular this implies a CLT\ for functions of reversible Markov
chains under the normalization $\sqrt{n}$. Kipnis and Varadhan (1986) provided
a functional form of this result. They showed that for a stationary reversible
and ergodic Markov chain the condition $var(S_{n})/n\rightarrow\sigma_{f}^{2}$
implies the convergence of $S_{[ns]}/\sqrt{n}$ to the Brownian motion
$|\sigma_{f}|W(s),$ (here, $0\leq s\leq1,$ $[ns]$ is the integer part of $ns$
and $W(s)$ is the standard Brownian motion). All these results used the
normalization $\sqrt{n.}$

Zhao et al. (2010) addressed the central limit theorem question for reversible
Markov chains under the weaker condition,
\begin{equation}
\sigma_{n}^{2}=nh(n), \label{var}%
\end{equation}
where $h$ is a slowly varying function (i.e. $\lim_{n\rightarrow\infty
}h(nu)/h(n)=1$ for all $u>0$). First, they proved that if $\lim\inf
_{n}h(n)=0,$ then necessarily
\begin{equation}
2S_{n}=(1+(-1)^{n-1}){X_{1\ }\ a.s.} \label{zero}%
\end{equation}
and if $\lim\inf_{n}h(n)=c^{2}\neq0$ then $\sigma_{n}^{2}/n\rightarrow c^{2}.$
In the following, to avoid the trivial case, we shall assume that
\[
\lim\inf\sigma_{n}^{2}/n>0.
\]
Zhao et al. (2010) also showed, by a class of examples, the surprising result
that the distribution of $S_{n}/\sigma_{n}$ needs not converge to the standard
normal distribution under (\ref{var}). Their example satisfies the central
limit theorem, in the sense that
\begin{equation}
\frac{S_{n}}{\sigma_{n}}\Rightarrow|c|Z, \label{CLT}%
\end{equation}
where $Z$ is a standard normal variable,$\ N(0,1),$ and $|c|\neq0,1$. A large
class of examples satisfying (\ref{CLT}) is given in Deligiannidis et al.
(2014). This paper also contains necessary and sufficient conditions for
(\ref{var}) in terms of the operator spectral measure.

A natural question is whether, in the context of reversible Markov chains, the
central limit theorem in (\ref{CLT}) implies the invariance principle namely
\begin{equation}
\frac{S_{[ns]}}{\sigma_{n}}\Rightarrow|c|W(s), \label{IP}%
\end{equation}
where $W(s),$ $s\geq0$ is the standard Brownian motion. This question is
interesting in itself, especially in the light of recent examples of
stationary sequences which satisfy CLT but not its functional forms (Giraudo
and Voln\'{y}, 2014).

A step in this direction is Theorem 2 in Longla et al. (2012), showing that
conditional CLT implies the functional CLT in (\ref{IP}) for functions of
reversible Markov chains. By the conditional CLT we understand that some $c>0$
and for all $t,$%
\begin{equation}
\mathbb{E}_{\xi_{0}}(\exp it\frac{S_{n}}{\sigma_{n}})\rightarrow\exp
(-\frac{t^{2}}{2c^{2}})\text{ in probability.} \label{cond CLT}%
\end{equation}
This form of the conditional CLT\ was essentially used in Longla et al.
(2012), in order to establish the convergence of finite dimensional distributions.

In this paper we show that actually (\ref{CLT}) and (\ref{var})\textit{
}implies (\ref{IP}) for functionals of reversible Markov chains. More
precisely we shall establish the following theorem:

\begin{theorem}
\label{main}Assume that $(X_{n})$ is defined by (\ref{defX}) and the Markov
chain is stationary and reversible. Then, the following statements are equivalent:
\end{theorem}

\textit{(a)\ The functional CLT\ in (\ref{IP}) holds.}

\textit{(b) The CLT\ in (\ref{CLT}) holds and the variance of partial sums is
regularly varying with exponent }$1$\textit{ (as in relation (\ref{var})).}

\textit{(c) The conditional CLT\ in (\ref{cond CLT}) holds.}

\section{Proof of Theorem \ref{main}}

The fact that (a) implies (b) follows by standard arguments in the following
way. Clearly, since the partial sum is just a finite dimensional distribution
for $s=1$, the functional CLT in (\ref{IP})\ implies the CLT in\ (\ref{CLT}).
Then, by (\ref{CLT})\textit{,}\ we have that $S_{[ns]}/\sigma_{\lbrack
ns]}\Rightarrow|c|Z,$ for every $s>0.$ On the other hand, the convergence in
(\ref{IP}) implies $S_{[ns]}/\sigma_{n}\Rightarrow|c|sZ$ for all $s>0$. By the
theorem of types (see Theorem 14.2 in Billingsley, 1995), it follows that
$\sigma_{\lbrack ns]}^{2}/\sigma_{n}^{2}\rightarrow s.$ The fact that (c)
implies (a) was established in Theorem 2 in Longla et al. (2012). It remains
to prove that (b) implies (c). The idea of proof is to show that from any
subsequence of $\mathbb{E}_{\xi_{0}}(\exp it(S_{n}/\sigma_{n}))$ we can
extract one converging to $\exp(-t^{2}/2c^{2})$ in probability. To achieve
this goal, we need a technical lemma concerning conditional convergence.

\begin{lemma}
\label{condlemma}Assume $(V_{n},\eta)$ is convergent in distribution to
$(V,Y)$. Then, we can construct on the same probability space a sequence
$(V_{n}^{\prime},\eta^{\prime}),$ where each vector is distributed as
$(V_{n},\eta)$ and a vector $(V^{\prime},\eta^{\prime})$ distributed as
$(V,Y)$ such that for all $t,$
\[
\mathbb{\ E}_{\eta^{\prime}}\exp(itV_{n}^{\prime})\rightarrow\mathbb{E}%
_{\eta^{\prime}}\exp(itV^{\prime})\text{ in probability.}%
\]

\end{lemma}

\textbf{Proof of Lemma \ref{condlemma}}. By the Skorohod theorem (see
Skorohod, 1956), we can construct on the same probability space a sequence
$(V_{n}^{\prime},\eta^{\prime}),$ where each vector is distributed as
$(V_{n},\eta),$ and a vector $(V^{\prime},Y^{\prime})$ distributed as $(V,Y)$
such that $(V_{n}^{\prime},\eta^{\prime})(\omega)\rightarrow(V^{\prime
},Y^{\prime})(\omega)$ for all $\omega.$ Clearly $\eta^{\prime}=Y^{\prime}$
and $V_{n}^{\prime}(\omega)\rightarrow V^{\prime}(\omega)$ for all $\omega.$
Then, by the mean value theorem, for any $t$ and $\delta>0$,
\begin{gather*}
|\mathbb{E}_{\eta^{\prime}}\exp(itV_{n}^{\prime})-\mathbb{E}_{\eta^{\prime}%
}\exp(itV^{\prime})|\leq\\
\mathbb{E}_{\eta^{\prime}}|\exp(itV_{n}^{\prime})-\exp(itV^{\prime}%
)|I(|V_{n}^{\prime}-V^{\prime}|\leq\delta)+2\mathbb{P}_{\eta^{\prime}}%
(|V_{n}^{\prime}-V^{\prime}|>\delta)\\
\leq|t|\delta+2\mathbb{P}_{\eta^{\prime}}(|V_{n}^{\prime}-V^{\prime}%
|>\delta)\text{ a.s.}%
\end{gather*}
By taking the expectation%
\[
\mathbb{E}|\mathbb{E}_{\eta^{\prime}}\exp(itV_{n}^{\prime})-\mathbb{E}%
_{\eta^{\prime}}\exp(itV^{\prime})|\leq|t|\delta+2\mathbb{P}(|V_{n}^{\prime
}-V^{\prime}|>\delta),
\]
which tends to $0$ as $n\rightarrow\infty$ and then $\delta\rightarrow0$.
$\ \ \ \ \diamond$

\bigskip

We continue to prove Theorem \ref{main} by proving that (b) implies (c). We
remind that we work under the assumption that $\lim\inf$ $\sigma_{n}^{2}/n>0.$
Note that, by stationarity and the fact that $X_{0}$ is square integrable, it
follows that $\max_{1\leq i\leq n}|X_{i}|/\sqrt{n}\rightarrow0$ $a.s.$ and in
$\mathbb{L}^{2}$ and therefore%
\begin{equation}
\max_{1\leq i\leq n}|X_{i}|/\sigma_{n}\rightarrow0\text{ }a.s.\text{ and in
}\mathbb{L}^{2}. \label{negl}%
\end{equation}
This property will allow us to adjust the sums by a few variables without
changing the limiting distribution. We shall use some notations: $\bar{S}%
_{n}=X_{n+1}+...+X_{2n}$; $\mathcal{P}_{n}=\sigma(\xi_{i},i\leq n)$ is the
past sigma field.

\textbf{Step 1}. As a preliminary computation, we show that for all $t,$
\begin{equation}
\mathbb{E(E}_{\xi_{0}}(\exp\frac{itS_{n}}{\sigma_{n}}))^{2}\rightarrow
\exp(-\frac{t^{2}}{c^{2}})\text{ .} \label{claim}%
\end{equation}
By (\ref{negl}), the properties of conditional expectation and Markov
property,
\begin{gather*}
\mathbb{E}(\exp\frac{itS_{2n}}{\sigma_{n}})=\mathbb{E}(\exp\frac{it(S_{n}%
+\bar{S}_{n})}{\sigma_{n}})=\mathbb{E}\left(  (\exp\frac{itS_{n}}{\sigma_{n}%
})\mathbb{E}(\exp\frac{it\bar{S}_{n}}{\sigma_{n}}|\mathcal{P}_{n})\right) \\
=\mathbb{E}\left(  (\exp\frac{itS_{n}}{\sigma_{n}})\mathbb{E}_{\xi_{n}}%
(\exp\frac{it\bar{S}_{n}}{\sigma_{n}})\right)  .
\end{gather*}
We write now $S_{n}=(X_{0}+...+X_{n-1})-X_{0}+X_{n}=S_{[0,n-1]}+(X_{n}%
-X_{0}).$

Clearly by simple computations and (\ref{negl}),%
\begin{gather*}
{\Huge |}\mathbb{E}\left(  \exp\frac{itS_{n}}{\sigma_{n}}\mathbb{E}_{\xi_{n}%
}(\exp\frac{it\bar{S}_{n}}{\sigma_{n}})\right)  -\mathbb{E}\left(  \exp
\frac{itS_{[0,n-1]}}{\sigma_{n}}\mathbb{E}_{\xi_{n}}(\exp\frac{it\bar{S}_{n}%
}{\sigma_{n}})\right)  {\Huge |}\\
\leq\mathbb{E}|\exp\frac{it(X_{n}-X_{0})}{\sigma_{n}}-1|\rightarrow0\text{ as
}n\rightarrow\infty.
\end{gather*}
Therefore, by combining these facts we obtain
\begin{equation}
\mathbb{E}(\exp\frac{itS_{2n}}{\sigma_{n}})-\mathbb{E}\left(  \exp
\frac{itS_{[0,n-1]}}{\sigma_{n}}\mathbb{E}_{\xi_{n}}(\exp\frac{it\bar{S}_{n}%
}{\sigma_{n}})\right)  \rightarrow0. \label{conv}%
\end{equation}
Note now that, by the properties of conditional expectation,%
\[
\mathbb{E}\left(  \exp\frac{itS_{[0,n-1]}}{\sigma_{n}}\mathbb{E}_{\xi_{n}%
}(\exp\frac{it\bar{S}_{n}}{\sigma_{n}})\right)  =\mathbb{E}\left(
\mathbb{E}_{\xi_{n}}(\exp\frac{itS_{[0,n-1]}}{\sigma_{n}})\mathbb{E}_{\xi_{n}%
}(\exp\frac{it\bar{S}_{n}}{\sigma_{n}})\right)  .
\]
By taking into account the reversibility of the process, we obtain
\[
\mathbb{E}_{\xi_{n}}(\exp\frac{itS_{[0,n-1]}}{\sigma_{n}})=\mathbb{E}_{\xi
_{n}}(\exp\frac{it\bar{S}_{n}}{\sigma_{n}}).
\]
Therefore, by combining this identity with (\ref{conv}), it follows that%
\[
\mathbb{E}(\exp\frac{itS_{2n}}{\sigma_{n}})-\mathbb{E}\left(  \mathbb{E}%
_{\xi_{n}}(\exp\frac{it\bar{S}_{n}}{\sigma_{n}})\right)  ^{2}\rightarrow0.
\]
So, by using now the stationary, the representation of the variance in
(\ref{var}) and the central limit theorem in (\ref{CLT}), the convergence
given in relation (\ref{claim}) follows.

\textbf{Step 2}{. }We show now that from any subsequence $(n^{\prime})$ we can
extract one $(n")\subset(n^{\prime})$ for which there is a pair of random
variables $(V,\xi),$ with $\xi$ distributed as $\xi_{0}$ and $V$ is normally
distributed with mean $0$ and variance $c^{2}$ (i.e. $N(0,c^{2})$) and such
that
\begin{equation}
\mathbb{\ E}_{\xi}\exp(itV_{n"}^{\prime})\rightarrow\mathbb{E}_{\xi}%
\exp(itV)\text{ in probability,} \label{subseq}%
\end{equation}
where $(V_{n"}^{\prime},\xi)$ is distributed as $(S_{n"}/\sigma_{n"},\xi
_{0}).$ In addition,
\begin{equation}
\mathbb{E(E}_{\xi_{0}}(\exp\frac{itS_{n"}}{\sigma_{n"}}))^{2}\rightarrow
\mathbb{E(E}_{\xi}\exp(itV))^{2}. \label{subsecE}%
\end{equation}
Indeed, since $(S_{n}/\sigma_{n},\xi_{0})$ is tight, from any subsequence
$(n^{\prime})$ we can extract one $(n")\subset(n^{\prime})$ such that
$(S_{n"}/\sigma_{n"},\xi_{0})$ is convergent in distribution to $(V,Y)$ where,
by (\ref{CLT}),\ $V$ is $N(0,c^{2})$ and $Y$ is distributed as $\xi_{0}$. By
Lemma \ref{condlemma}, applied to $V_{n"}=S_{n"}/\sigma_{n"},$ there exists
two pairs of variables: $(V_{n"}^{\prime},\xi)$ distributed as $(S_{n"}%
/\sigma_{n"},\xi_{0}^{\prime})$ and $(V,\xi),$ with $\xi\ $distributed as
$\xi_{0}$ and $V$ centered normal with variance $c^{2}$, such that convergence
in (\ref{subseq})\ holds. Now by starting from (\ref{subseq}) and applying the
Lebesgue dominated convergence theorem, (\ref{subsecE}) follows.

\textbf{Step 3}. In order to finish the proof of the theorem we shall show
that the limit in (\ref{subseq}) does not depend on the subsequence. As a
matter of fact we shall show that
\begin{equation}
\mathbb{E}_{\xi}\exp(itV)=\exp(-\frac{t^{2}}{2c^{2}})\text{ a.s.} \label{CF}%
\end{equation}
We shall use the fact that $V$ is $N(0,c^{2})$ and also, by Step $2$ and by
Step $1$, we know that
\begin{equation}
\mathbb{E(E}_{\xi}\exp(itV))^{2}=\exp(-\frac{t^{2}}{c^{2}})\text{.}
\label{normalid}%
\end{equation}
Note that, in order to establish (\ref{CF}), it is enough to show that, for
all integers $m\geq0,$
\begin{equation}
\mathbb{E}_{\xi}(V^{m})=\mathbb{E}(V^{m})\;a.s. \label{H}%
\end{equation}
With this aim, we redefine $V$ on a larger probability space together with two
independent variables uniformly distributed $U$ and $\tilde{U}$ which are also
independent on $\xi$ and such that $V=f(\xi,U)$ and $\tilde{V}=f(\xi,\tilde
{U}).$ Note that $V$ and $\tilde{V}$ are conditionally independent given
$\xi\ $and $(V,\xi)$ has the same distribution as $(\tilde{V},\xi).$ In
addition, let $N,\tilde{N}$ be i.i.d. random variables $N(0,c^{2})$.

Then, by (\ref{normalid}),
\begin{gather*}
\mathbb{E(E}_{\xi}\exp(itV))^{2}=\mathbb{E(}\exp(it(V+\tilde{V})))^{2}\\
=\sum_{m=0}^{\infty}\mathbb{E}(V+\tilde{V})^{m}\frac{(it)^{m}}{m!}=\sum
_{m=0}^{\infty}\mathbb{E}(N+\tilde{N})^{m}\frac{(it)^{m}}{m!}.
\end{gather*}
Hence, for all $n=0,1,2,\ldots$
\begin{equation}
\mathbb{E}(V+\tilde{V})^{n}=\mathbb{E}(N+\tilde{N})^{n}. \label{ID}%
\end{equation}
Further, we proceed by induction to prove (\ref{H}). Note that (\ref{H})
obviously holds for $m=0$. Assume (\ref{H}) holds for $m\leq k$. To prove it
for $(k+1)$ we use (\ref{ID}) with $n=2k+2$ and develop the binomial%

\begin{equation}
\sum_{\ell=0}^{2k+2}C_{2k+2}^{\ell}\mathbb{E}V^{\ell}\tilde{V}^{2k+2-\ell
}=\sum_{\ell=0}^{2k+2}C_{2k+2}^{\ell}\mathbb{E}N^{\ell}\mathbb{E}N^{2k+2-\ell
}. \label{binom}%
\end{equation}
By the induction hypothesis, conditional independence of $V$ and $\tilde{V}$
and the properties of conditional expectation we obtain, for $0\leq\ell\leq
k,$ and every integer $m,$%
\[
\mathbb{E}V^{\ell}\tilde{V}^{m}=\mathbb{E(E}_{\xi}(V^{\ell}\tilde{V}%
^{m}))=\mathbb{E(E}_{\xi}V^{\ell}\mathbb{E}_{\xi}V^{m})=\mathbb{E(}V^{\ell
})\mathbb{E(}V^{m})=\mathbb{E(}N^{\ell})\mathbb{E(}N^{m})
\]
and similarly,
\[
\mathbb{E}\tilde{V}^{\ell}V^{m}=\mathbb{E(}N^{\ell})\mathbb{E(}N^{m}).
\]
By using both these estimates in (\ref{binom}), we observe that the terms in
the sum from $0\leq\ell\leq k$ and $k+2\leq\ell\leq2k+2\ $all cancel, and it
follows that
\[
\mathbb{E}V^{k+1}\tilde{V}^{k+1}=(\mathbb{E}N^{k+1})^{2}=(\mathbb{E}%
V^{k+1})^{2}.
\]
Taking into account that
\[
\mathbb{E}V^{k+1}\tilde{V}^{k+1}=\mathbb{E(E}_{\xi}V^{k+1}\mathbb{E}_{\xi
}\tilde{V}^{k+1})=\mathbb{E(E}_{\xi}V^{k+1})^{2},
\]
we obtain by the above arguments that
\[
\mathbb{E(E}_{\xi}V^{k+1})^{2}=(\mathbb{E}V^{k+1})^{2}.
\]
It follows that%
\[
\mathbb{E(E}_{\xi}V^{k+1}-\mathbb{E}V^{k+1})^{2}=0,
\]
which completes the proof of (\ref{H}) and therefore of (\ref{CF}). Combining
(\ref{CF}) with (\ref{subseq}) we get $\mathbb{\ }$%
\[
\mathbb{E}_{\xi}\exp(\frac{itS_{n"}^{\prime}}{\sigma_{n"}})\rightarrow
\exp(-\frac{t^{2}}{2c^{2}})\text{ in probability.}%
\]
Consequently,
\[
\mathbb{E}_{\xi_{0}}\exp(\frac{itS_{n"}}{\sigma_{n"}})\rightarrow\exp
(-\frac{t^{2}}{2c^{2}})\text{ in probability,}%
\]
completing the proof of the theorem.

\bigskip

\textbf{Acknowledgement}. This \ material was supported in part by a Charles
Phelps Taft Memorial Fund grant at the University of Cincinnati, and the
National Science Foundation grant number DMS-1208237. The authors would like
to thank the referee for carefully reading the manuscript and suggestions
which improved the presentation of the paper.

\end{document}